\documentstyle[amssymb,amstex,graphicx,graphics,12pt,a4wide]{article}
\newtheorem{thm}{Theorem}
\newtheorem{lemma}{Lemma}
\newtheorem{proposition}{Proposition}

\newtheorem{remark}{Remark}

\newcommand{\Vol}{\operatorname{Vol}}

\begin{document}

\title{\textbf{Exact and asymptotic results for intrinsic volumes of Poisson $k$-flat processes}}
\author{Matthias Schulte\footnote{\textit{e-mail: matthias.schulte[at]uni-osnabrueck.de}} and Christoph Th\"ale\footnote{\textit{e-mail: christoph.thaele[at]uni-osnabrueck.de}}\\
Department of Mathematics\\
University of Osnabr\"uck, Germany}
\date{}
\maketitle

\begin{abstract}
The intrinsic volumes induced by a stationary Poisson $k$-flat process inside a compact and convex sampling window are considered. Using techniques from stochastic analysis, more precisely calculus with multiple stochastic integrals and a Wiener-It\^o chaos expansion of Poisson functionals, all moments and cumulants, exact and asymptotic, are calculated in terms of a family of integral-geometric functionals. Moreover, univariate central limit theorems as well as Berry-Esseen-type inequalities are shown and a multidimensional limit theorem is concluded.
\end{abstract}
\begin{flushleft}\footnotesize
\end{flushleft}

\section{Introduction}

Random subdivisions of the $d$-dimensional Euclidean space (with $d\geq 2$) induced by a collection of random hyperplanes are one of the most classical and well studied models for random tessellations considered in stochastic geometry and spatial statistics. A class that has attracted particular interest in the literature is that of Poisson hyperplanes, which was introduced by Miles \cite{Miles0} and Matheron \cite{Materon2} and was studied by many authors (see \cite{SW} and \cite{SKM} and the references therein). Also the more general class of $k$-flats ($0\leq k\leq d-1$), in particular Poisson $k$-flats, was considered in the literature \cite{BL, HLW, SW, SP1, SP2}. However, only a very few explicit and non-asymptotic results are available for these models. Main contributions in this direction are due to Heinrich et al. \cite{HSS} and Heinrich \cite{H} in the hyperplane case. The aim of the present paper is to exploit a new technique resting upon a Wiener-It\^o-chaos decomposition and Fock space representation of Poisson functionals \cite{LP} and to show how it can successfully be applied to get new and deep insight into the structure of Poisson $k$-flats in ${\Bbb R}^d$. We are able to extend the known results from the literature in two directions. On the one hand and in contrast to \cite{H, HSS} we will not only deal with hyperplanes, but with more general $k$-flat processes and consider certain lower-order intrinsic volumes, whereas \cite{H, HSS} only deal with $k$-volume and Euler-characteristic, two special intrinsic volumes, which are included in our discussion. On the other hand, we can derive elegant \textit{exact} and \textit{non-asymptotic} formulae for all moments and cumulants -- completely characterizing thereby also distributions -- in terms of a class of integral-geometric functionals. The precise statements of our results are the content of the next section. %The functionals considered in this paper have a in some sense simple form, making the resulting exact expressions particularly elegant. In our further work in progress we will demonstrate that the results get much more involved when we consider lower-dimensional intersection processes.
\\ 
As already signalled above, our results rest upon a Wiener-It\^o chaos expansion of Poisson functionals. This is a very powerful and advanced technique from stochastic analysis that has not been used so far for problems in stochastic geometry. To keep the paper self-contained and for easier reference we review in Section \ref{secCHAOS} some of the most important results related to that theory, which are needed for our proofs. They are the content of Section \ref{secPROOFS} and the key to them is a product formula for multiple stochastic integrals. Such a formula is available in the literature in great generality, see \cite{PT, Sur}. However, for our purposes a special and easily understandable particular case is sufficient. For this reason and in order to keep the argument simple, elegant and transparent, we circumvent the use of such general results and derive a recursion formula for the product of a special class of multiple Wiener-It\^o integrals with respect to Poisson processes, which is much easier to handle than the classical product or diagram formula.\\ We would like to point out that such an elegant and non-technical way is no more possible when dealing with lower-dimensional intersection processes and their induced intrinsic volumes. Their study is for this reason postponed to a future separate paper.

\section{Statement of results}\label{secIntro}

In this paper, we study stationary Poisson $k$-flat processes $\eta_k$ in ${\Bbb R}^d$ ($0\leq k\leq d-1$) with intensity measure $\tau_k\Lambda_k$, where $\tau_k\in(0,\infty)$, that is to say $\eta_k$ is a stochastically translation invariant point process on the space ${\cal E}_k$ of $k$-flats (this is the space of $k$-dimensional affine planes) in ${\Bbb R}^d$ and the mean $k$-volume ${\Bbb  E}\mbox{Vol}_{k}(\eta_k\cap T)$ (here and in the sequel $\Vol_k$ stands for the $k$-dimensional volume measure) generated in a test set $T\subset{\Bbb  R}^d$ of unit volume equals $\tau_k$ (as usual, we identify a point process with the random closed set induced by it). Here and in the sequel, $\Lambda_k$ stands for a unit density translation invariant measure on ${\cal E}_k$, by which we mean a measure of the product form $\ell\otimes{\cal R}_k$, where $\ell$ stands for the Lebesgue measure on ${\Bbb R}$ and ${\cal R}_k$ for a probability measure on the space ${\cal L}_k$ of $k$-dimensional linear subspaces of ${\Bbb R}^d$. If ${\cal R}_k$ is the uniform distribution on ${\cal L}_k$ with probability element $dL_k$, then $\eta_k$ corresponds to a stationary and isotropic (this is stochastically rotation invariant) Poisson $k$-flat process. The focus of our considerations is on the intrinsic volumes $V_{j,k}(W)=V_j(\eta_k\cap W)$, $0\leq j\leq k$, $\eta_k$ induces in a $d$-dimensional compact convex body $W\subset{\Bbb  R}^d$ with interior points. For a convex set $K\subset{\Bbb R}^d$ they are defined by $$V_d(K)=\Vol_d(K),\ \ \ V_j(K)={{d\choose j}\kappa_d\over\kappa_j\kappa_{d-j}}\int_{{\cal L}_j}\Vol_j(K|L_j)dL_j,\ \ \ \ 0\leq j\leq d-1$$ where $K|L_j$ stands for the orthogonal projection of $K$ onto $L_j$ and where $\kappa_i$ stands for the volume of the $i$-dimensional unit ball. If the affine hull of $K\subset{\Bbb R}^d$ is $k$-dimensional, $V_k(K)$ is the $k$-volume, $2V_{k-1}(K)$ the $(k-1)$-dimensional surface area, ${2\kappa_{k-1}\over k\kappa_k}V_1(K)$ the mean width and $V_0(K)=1$ the Euler-characteristic of $K$. By additive extension, the intrinsic volumes are also well defined for locally finite unions of convex sets, cf. \cite{SW}. In this sense $V_{k,k}(W)$ is just the $k$-volume or $V_{0,k}(W)$ the number of $k$-flats of $\eta_k$ in $W$. We will also consider a sequence $W_\varrho=\varrho W$ of dilated windows in which $\eta_k$ is observed.\\
%and can hence assume without loss of generality $W\subset B^d_r$ with $r=\kappa_{d-1}^{-\frac{1}{d-1}}$, where $\kappa_{d-1}$ is the volume of the $(d-1)$-dimensional unit ball.
%The quantity $V_{j,k}(W)$ is of interest for example in material sciences. More precisely, imagine material modeled by a system of Poisson hyperplanes (here $k=d-1$), then $V_{j,d-1}(W)$ encodes important properties of this material. For example, the cell boundaries could be semi-permeable and a fluid could perfuse the material. In this case, the knowledge of the total surface content of the hyperplane system is important for the analysis of the velocity of the flow or the total fluid volume passing through the material. Other applications related to recent developments in telecommunication were discussed for example in \cite{HSS}.\\
Under our invariance assumption, it is readily seen that ${\Bbb E}V_{k,k}(W)=\tau_k\mbox{Vol}_d(W)$. In sharp contrast to the mean values, higher-order moments will depend on the precise shape of the window $W$ via integral-geometric functionals $A(W,j,k,m)$ defined by $$A(W,j,k,m):=\int_{[W]_k}V_j^m(W\cap E_k)\Lambda_k(dE_k),\ \ \ \ m=1,2,3,\ldots,\ \ 0\leq j\leq k,$$ where $[W]_k\subset{\cal E}_k$ denotes the set of $k$-flats hitting $W$. Let us denote by $\mu_{m}(V_{j,k}(W))$ the $m$-th moment of $V_{j,k}(W)-{\Bbb  E}V_{j,k}(W)$, i.e. $\mu_m(V_{j,k}(W)):={\Bbb  E}(V_{j,k}(W)-{\Bbb  E}V_{j,k}(W))^m$ and by $\gamma_m(V_{j,k}(W))$ the cumulant of order $m$ (recall that the cumulants of a random variable are the coefficients of the series expansion of the logarithm of its moment generating function, see e.g. \cite[p. 290]{Sir} or (\ref{cumulant}) below), in particular ${\Bbb V}V_{j,k}(W)=\gamma_2(V_{j,k}(W))$. We also consider the normalized intrinsic volumes $V_{j,k}^*(W)$ defined by $$V_{j,k}^*(W):={V_{j,k}(W)-{\Bbb E}V_{j,k}(W)\over\sqrt{{\Bbb V}V_{j,k}(W)}}.$$ To neatly formulate our results, denote by $\Pi(\left\{1,\ldots,m\right\})$ the set of partitions of $\left\{1,\ldots,m\right\}$ in disjoint sets with at least two elements and by $\Pi_m$ the multiset $$\Pi_m=\left\{\left\{m_1,\ldots,m_r\right\}:(A_i)_{i=1}^r\in\Pi(\left\{1,\ldots,m\right\}),m_i=|A_i|\mbox{ for }i=1,\ldots,r\right\}.$$ The cardinality of $\left\{m_1,\ldots,m_r\right\}$ in the  multiset $\Pi_m$ is given by the multinomial coefficient ${m \choose m_1 \ldots m_r}$.
\begin{thm}\label{thm1} Consider a stationary Poisson $k$-flat process in ${\Bbb R}^d$, $0\leq k\leq d-1$ as above.
\begin{itemize}
 \item[(a)] For a fixed compact convex body $W\subset{\Bbb R}^d$ with interior points we have $\mu_1(V_{j,k}(W))=\gamma_1(V_{j,k}(W))=0$, $$\mu_m(V_{j,k}(W))=\sum_{\pi\in\Pi_m}\tau_k^{|\pi|}\prod_{m_l\in\pi}A(W,j,k,m_l)$$ and $$\gamma_m(V_{j,k}(W))=\tau_k A(W,j,k,m)$$ for $m\geq 2$ and explicit expressions for $A(W,j,k,m)$ are provided in Lemma \ref{lemmaABdn} below for the case of a Euclidean ball.
 \item[(b)] For the growing sequence $W_\varrho=\varrho W$ the following holds: $$\lim_{\varrho\rightarrow\infty}{\mu_m(V_{j,k}(W_\varrho))\over\varrho^{jm+\lfloor{m(d-k)\over 2}\rfloor}}=\begin{cases}(m-1)!! \left(\tau_k A(W,j,k,2)\right)^{\frac{m}{2}} &: m\ \mbox{even}\\ {m\choose 3}(m-4)!!\tau_k A(W,j,k,3)(\tau_k A(W,j,k,2))^{\frac{m-3}{2}}&: m\ \mbox{odd}\end{cases}$$ for $m\geq 2$ with the double factorial $k!!=k(k-2)(k-4)\ldots$ and $$\lim_{\varrho\rightarrow\infty}{\gamma_m(V_{j,k}(W_\varrho))\over\varrho^{jm+d-k}}=\tau_k A(W,j,k,m).$$
 \item[(c)] For the normalized intrinsic volumes $V_{j,k}^*(W_\varrho)$ we have $$\lim_{\varrho\rightarrow\infty}{\mu_m(V_{j,k}^*(W_\varrho))}=\begin{cases}(m-1)!! &: m\ \mbox{even}\\ 0 &: m\ \mbox{odd}\end{cases}$$ for $m\geq 2$ and $$\lim_{\varrho\rightarrow\infty}\gamma_m(V_{j,k}^*(W_\varrho))=\begin{cases}1 &: m=2\\ 0 &: m=1\ \mbox{or}\ m\geq 3.\end{cases}$$
\end{itemize}
\end{thm}
%\begin{remark}
%Formulas for the moments and cumulants $\mu_m(V_{j,k}(W))$ and $\gamma_m(V_{j,k}(W))$ could also be concluded from Corollary 7.2 in \cite{PT}. However, this result is stated only for simple functions there in order to avoid integrability problems, which do not occur in our setting. Its proof is based on a complicated and rather technical product formula for a broad class of stochastic integrals, which was introduced in \cite{Sur} for the Poisson case. In contrast, we prefer to give an elementary and direct proof of Theorem \ref{thm1} (a), which only makes use of simple properties of multiple Wiener-It\^o integrals.
%\end{remark}
In general, the integral-geometric functionals $A(W,j,k,m)$ are very difficult to evaluate explicitly. However, for a stationary and isotropic Poisson $k$-flat process $\eta_k$ and the $d$-dimensional ball $B_\varrho^d$ with radius $\varrho>0$ we have the following closed formula:
\begin{lemma}\label{lemmaABdn} For an isotropic Poisson $k$-flat process, $W=B_\varrho^d$, $0\leq j\leq k\leq d-1$ and $m=1,2,3,\ldots$ we have
\begin{equation}\nonumber A(B_\varrho^d,j,k,m)=\left[{k\choose j}{\kappa_k\over\kappa_{k-j}}\right]^m{\kappa_{km+1}\over\kappa_{km}}\varrho^{jm+d-k}.\label{eqABdn}
\end{equation}
\end{lemma}
Thus, in the isotropic case and for $W=B^d:=B_1^d$, the moments and cumulants of $V_{j,k}(W)-{\Bbb  E}V_{j,k}(W)$ can explicitly be evaluated. It is also worth pointing out that in the isotropic set-up the special functional $A(W,j,k,2)$ may be interpreted as the order $j+1$ chord-power integral of $W$. Indeed, a repeated application of the affine Blaschke-Petkantschin formula from integral geometry shows the equality $$A(W,j,k,2)={\kappa_k\over j+1}\int_{{\cal L}_1}\mbox{Vol}_1^{j+1}(W\cap L_1)dL_1,$$ where ${\cal L}_1$ is the space of lines in ${\Bbb  R}^d$ with invariant measure $dL$, see \cite[Eq. (8.57)]{SW}.\\ It is easily verified that the moments $\mu_m(V_{j,k}(W))$ satisfy the assumption of the well known Carleman criterion for the uniqueness of the moment problem \cite[p. 296]{Sir}. Thus, the entire distribution of the centred intrinsic volumes -- and hence the intrinsic volumes itself -- is determined by the sequence $\mu_m(V_{j,k}(W))$. Further, the asymptotics for the moments show how the leading term of their expansion is influenced by the shape of the window $W$.\\ The knowledge of the distribution of $V_{j,k}(W)$ allows alternative proofs of a central limit theorem for $V_{j,k}^*(W_\varrho)$, a topic which was previously discussed in \cite{H, HSS, Paroux} in the special cases $k=d-1$ and $j=0$ or $j=k$. It is readily seen that the well known method of moments can be applied and also the celebrated cumulant method. Alternatively, the general estimate for the Wasserstein distance recently obtained in \cite{P} by means of Stein's method combined with Malliavin calculus for Poisson point processes can be used to conclude the central limit theorem. Moreover, we are in the position to prove a multivariate central limit theorem for the vector $(V_{0,k}^*(W_\varrho),\ldots,V_{k,k}^*(W_\varrho))$ of intrinsic volumes of a stationary Poisson $k$-flat process, which is new.
\begin{thm}\label{thm2} Let $\eta_k$ be a stationary $k$-flat process in ${\Bbb R}^d$, $0\leq k\leq d-1$.
\begin{itemize}
 \item[(a)] As $\varrho\rightarrow\infty$, $V_{j,k}^*(W_\varrho)$ converges in distribution to a standard normal random variable $\cal N$ and, moreover, $$\sup_{t\in{\Bbb R}}|{\Bbb P}(V_{j,k}^*(W_\varrho)\leq t)-{\Bbb P}({\cal N}\leq t)|\leq 42{A(W,j,k,2)\over A(W,j,k,3)}\tau_k^{-{1\over 2}}\varrho^{-{d-k\over 2}}.$$
 \item[(b)] The vector $(V_{0,k}^*(W_\varrho),\ldots,V_{k,k}^*(W_\varrho))$ converges in distribution to a $(k+1)$-dimensional Gaussian random vector with mean zero and covariance matrix $C=(c_{i,j})_{i,j=0}^k$ with $c_{i,j}$ given by $$c_{i,j}={\int_{[W]_k}V_i(W\cap E_k)V_j(W\cap E_k)\Lambda_k(dE_k)\over\sqrt{\int_{[W]_k}V_i^2(W\cap E_k)\Lambda_k(dE_k)\int_{[W]_k}V_j^2(W\cap E_k)\Lambda_k(dE_k)}}.$$ In the particular case $W=B^d$ we have $$c_{i,j}={\Gamma\left(1+{i+j\over 2}\right)\over\Gamma\left({3\over 2}+{i+j\over 2}\right)}\sqrt{{\Gamma\left({3\over 2}+i\right)\Gamma\left({3\over 2}+j\right)\over i!j!}}.$$
% \item[(b)] With $c=324\sqrt{2}$ it holds $$\sup_{x\in{\Bbb  R}}\left|{\Bbb  P}(S_*(W_\varrho)\leq x)-\Phi(x))\right|\leq{c\over \sqrt{\tau\varrho} A(W,2)^{3\over 2}},$$ where $\Phi$ is the distribution function of a standard normal random variable.
% \item[(c)] There exist constants $c_1,c_2>0$ such that for $x\in[0,c_1\sqrt{\tau\rho} A(W,2)^{3/2}]$ we have
%\begin{eqnarray}
%\nonumber \left|\log{{\Bbb  P}(S_*(W_\varrho)\geq x)\over 1-\Phi(x)}\right| &\leq & c_2{1+x^3\over\sqrt{\tau\varrho} A(W,2)^{3\over 2}},\\
%\nonumber \left|\log{{\Bbb  P}(S_*(W_\varrho)\leq -x)\over\Phi(-x)}\right| &\leq & c_2{1+x^3\over\sqrt{\tau\varrho} A(W,2)^{3\over 2}}.
%\end{eqnarray}
% \item[(d)] For all $x\geq 0$ we have $${\Bbb  P}(S_*(W_\varrho)\geq x)\leq\exp\left\lbrace -{1\over 4}\min\left\lbrace {x^2\over A(W,m)},x\sqrt{\tau\varrho} A(W,2)^{3\over 2}\right\rbrace\right\rbrace.$$
\end{itemize}
\end{thm}
\begin{remark}
Beside the multivariate normal convergence of part (b) in the previous theorem, also a rate of convergence for certain multivariate probability metrics can be concluded from \cite{PZ}. However, their definitions are quite involved and the corresponding results are omitted for this reason.
\end{remark}
%Theorem \ref{thm2} (a) not only ensures that $V_{j,k}^*(W_\varrho)$ asymptotically tends to a standard normal random variable $\cal N$, but also makes a precise statement about the rate of convergence. Moreover, this rate is geometry-sensitive in a sense that it includes information on how the shape of the window $W$ has influence on this rate.

\section{Wiener-It\^o Chaos Expansions}\label{secCHAOS}

This section is devoted to a brief introduction to Wiener-It\^o chaos expansions.  Let $({\cal X},{\cal B}({\cal X}),\lambda)$ be a standard Borel space, which means that there is a bijection between ${\cal X}$ and a Borel subset of $[0,1]$, which is measurable in both directions, see \cite{Kal}. Moreover, we assume that the measure $\lambda$ is $\sigma$-finite and has no atoms. Let $\eta$ be a Poisson point process with intensity measure $\lambda$ on an underlying probability space $(\Omega,{\cal F},{\Bbb  P})$ . By ${\Bbb  P}_\eta$ we denote the distribution of $\eta$ on the space $N({\cal X})$ of all integer-valued $\sigma$-finite measures on ${\cal X}$. The space $N({\cal X})$ is equipped with the smallest $\sigma$-field such that all maps $\nu\mapsto \nu(A)$ with $A\in{\cal B}({\cal X})$ and $\nu\in N({\cal X})$ are measurable.\\
We present now a definition for multiple Wiener-It\^o integrals, following thereby \cite{Sur}. One starts with simple functions and extends the definition to arbitrary functions in $L_s^2({\cal X}^n)$, $n\in{\Bbb N}$, the space of square-integrable symmetric functions with norm $||\cdot||_n$ and inner product $\langle\cdot,\cdot\rangle_n$. A function $f\in L^2({\cal X}^n)$ is called simple, if
\begin{enumerate}
 \item $f$ is symmetric,
 \item $f$ is constant on a finite number of Cartesian products $B_1\times\ldots\times B_n$ and vanishes elsewhere,
 \item $f$ vanishes on diagonals, that means $f(x_1,\ldots,x_n)=0$ if $x_m= x_l$ for some $m\neq l$.
\end{enumerate}
Let $L^2_0({\cal X}^n)$ be the space of all simple functions. For $f_0\in L^2_0({\cal X}^n)$, the multiple Wiener-It\^o integral $I_n(f_0)$ of $f_0$ with respect to the compensated Poisson point process $\eta-\lambda$ is defined as
$$I_n(f_0)=\int f_0d(\eta-\lambda)=\sum f_0^{B_1\times\ldots\times B_n}(\eta-\lambda)(B_1)\ldots(\eta-\lambda)(B_n),$$
where we sum over all Cartesian products and $f_0^{B_1\times\ldots\times B_n}$ is the constant value on such a set. A straight forward computation shows that
\begin{equation}{\Bbb  E}I_n(f_0)^2=n!||f_0||_n^2.\label{eqISO}\end{equation}
Thus, there is an isometry between $L^2_0({\cal X}^n)$ and a subset of $L^2({\Bbb  P}_\eta)$. Furthermore, $L^2_0({\cal X}^n)$ is dense in $L^2_s({\cal X}^n)$, whence for every $f\in L^2_s({\cal X}^n)$ there is a sequence $(f_l)_{l\in{\Bbb N}}$ of simple functions with $f_l\rightarrow f$ in $L^2_s({\cal X}^n)$. Because of the isometry (\ref{eqISO}), it is possible to define $I_n(f)$ as the limit of $(I_n(f_l))_{l\in{\Bbb N}}$ in $L^2({\Bbb  P}_\eta)$.\\
\begin{remark} The fact that $L_0^2({\cal X}^n)$ is dense in $L_s^2({\cal X}^n)$ depends on the assumed topological structure of the space ${\cal X}$ and the fact that $\lambda$ is non-atomic. For a definition without these requirements we refer to \cite{LP}.
\end{remark}
It follows directly from the definition that multiple Wiener-It\^o integrals have the properties summarized in the following
\begin{lemma}\label{lemItoIntegrals}
Let $f\in L^2_s({\cal X}^n)$ and $g\in L^2_s({\cal X}^m)$ with $m,n\geq 1$ and $n\neq m$. Then
$$\text{(a)}\ {\Bbb  E}I_n(f)=0,\ \ \ \ \ \text{(b)}\ {\Bbb  E}I_n(f)^2=n!||f||_n^2,\ \ \ \ \ \text{(c)}\ {\Bbb  E}I_n(f)I_m(g)=0.$$
\end{lemma}
The second ingredient for Wiener-It\^o chaos expansions is the difference operator $D$. Let $F: N({\cal X})\mapsto\overline{{\Bbb R}}:={\Bbb R}\cup\left\{\pm\infty\right\}$ be measurable and $x\in {\cal X}$. In this case, the difference operator (or add-one-cost operator) $D_x$ is defined by $$D_xF(\eta):=F(\eta+\delta_x)-F(\eta),$$ where $\delta_x$ stands for the unit mass Dirac measure concentrated at $x$. It is well known that the difference operator obeys the product formula
\begin{equation}
D_z(FG)=(D_zF)G+F(D_zG)+(D_zF)(D_zG), \label{eqProductDifference}
\end{equation}
see Lemma 6.1 in \cite{Nua}. The iterated difference operators $D_{x_1,\ldots,x_n}^n$ with $x_1,\ldots,x_n\in {\cal X}$ can now be introduced by $$D^0\equiv 1,\ \ \ \ \ D^1:=D,\ \ \ \ \ D_{x_1,\ldots,x_n}^nF(\eta):=D_{x_1}^1D_{x_2,\ldots,x_n}^{n-1}F(\eta).$$ Let us further define the functions $T_nF: {\cal X}^n\mapsto\overline{{\Bbb R}}$ 
via \begin{equation}\nonumber\label{eqKERNELS}
 T_0F={\Bbb  E}F(\eta),\ \ \ \ \ T_nF(x_1,\ldots,x_n)={\Bbb  E}D_{x_1,\ldots.x_n}^nF(\eta),\ \ \ \ n\geq 1,
\end{equation}
which are called the kernels of $F$ in the following. Because of the symmetry of the iterated difference operators, the kernels are symmetric, too. It was shown in \cite[Thm. 1.3]{LP} that any $F\in L^2({\Bbb  P}_\eta)$ can be decomposed uniquely as a sum of random variables that are the multiple Wiener-It\^o integrals of these kernels:
\begin{proposition}\label{propCHAOS} For $F\in L^2({\Bbb  P}_\eta)$ we have $T_nF\in L_s^2({\cal X}^n),n\geq 1$, and \begin{equation}F(\eta)=\sum_{n=0}^\infty{1\over n!}I_n(T_nF)\label{eqCHAOS}\end{equation} and the series converges in $L^2({\Bbb  P})$. Moreover, the kernels $T_nF$ are the $\lambda$-almost everywhere unique functions $g_n\in L^2_s({\cal X}^n)$, $n\in{\Bbb N}$, satisfying $F(\eta)=\sum_{n=1}^\infty \frac{1}{n!}I_n(g_n)$. Furthermore, \begin{equation}{\Bbb C}ov(F,G)=\sum_{n=1}^\infty{1\over n!}\langle T_nF,T_nG\rangle_n\label{eqn:CovChaos}\end{equation} for $F,G\in L^2({\Bbb P}_\eta)$.
\end{proposition}
The space of sequences $(f_n)_{n\in{\Bbb N}}$ with $f_n\in L^2_s({\cal X}^n)$ and $\sum_{n=1}^\infty \frac{1}{n!}||f_n||_n^2<\infty$ is often called the Fock space in the literature. Note that by Proposition \ref{propCHAOS}, we have an isometry between $L^2(P_\eta)$ and the Fock space. The identity (\ref{eqCHAOS}) is called Wiener-It\^o chaos expansion or Fock space representation of $F$.%\\ \\

\section{Proofs}\label{secPROOFS}

Recall that by $\eta_k$ we denote a stationary Poisson $k$-flat process with $k$-volume density $\tau_k$ and that for a $d$-dimensional compact convex body $W$, $V_{j,k}(W)=V_j(\eta_k\cap W)$ stands for the $j$-th intrinsic volume induced by $\eta_k$ in $W$. Moreover, $\mu_m(V_{j,k}(W))$ ist the $m$-th centred moment ${\Bbb  E}(V_{j,k}(W)-{\Bbb  E}V_{j,k}(W))^m$ and $\gamma_m(V_{j,k}(W))$ is the order $m$ cumulant of $V_{j,k}(W)$, i.e. \begin{equation}\log{\Bbb  E}e^{itV_{j,k}(W)}=\sum_{m=0}^\infty{(it)^m\over m!}\gamma_m(V_{j,k}(W)),\ \ \ \ t\in{\Bbb  R}.\label{cumulant}\end{equation} By $V_{j,k}^*(W)$ we denote the centred and normalized version $(V_{j,k}(W)-{\Bbb  E}V_{j,k}(W))/\sqrt{{\Bbb  V}V_{j,k}(W)}$. Recall further that by $W_\varrho$ we mean the family of dilated and growing windows $\varrho W\subset{\Bbb  R}^d$, $\varrho>0$.\\
In our exact moment formula, a crucial r\^ole is played by products of functionals $A(W,j,k,m)$. To derive an asymptotic expression it will be necessary to identify the leading $\varrho$-term. To do so, the following observation will be helpful.
\begin{proposition}\label{propHOMOA} The functional $A(W,j,k,m)$ is homogeneous of degree $jm+d-k$, i.e. $$A(W_\varrho,j,k,m)=\varrho^{jm+d-k}A(W,j,k,m),\ \ \ \varrho>0.$$
\end{proposition}
\paragraph{Proof} Noting that $\Lambda_k([W_\varrho]_k)=\varrho^{d-k}\Lambda_k([W]_k)$ for translation invariant measures $\Lambda_k$ on ${\cal E}_k$ implies
\begin{eqnarray}
\nonumber & & A(W_\varrho,j,k,m)\\
\nonumber &=& \int_{[W_\varrho]_k}V_j^m(W\cap E_k)\Lambda_k(dE_k) = \int_{[W]_k}V_j^m(\varrho(W\cap E_k))\varrho^{d-k}\Lambda_k(dE_k)\\
\nonumber &=& \varrho^{jm+d-k}\int_{[W]_k}V_j^m(W\cap E_k)\Lambda_k(dE_k) = \varrho^{jm+d-k}A(W,j,k,m),
\end{eqnarray}
where we have used that $V_j$ is homogeneous of degree $j$.\hfill $\Box$

\paragraph{\textsc{Proof of Lemma \ref{lemmaABdn}}} Recall that we take for $\Lambda_k$ the translation invariant measure on ${\cal E}_k$ induced by the unique unit density Haar measure on ${\cal L}_k$. In view of Proposition \ref{propHOMOA} it is sufficient to consider the case $\varrho=1$, where we write $B^d$ instead of $B_1^d$. Observe that the intersection $B^d\cap E_k(p)$ is a $k$-dimensional ball with radius $\sqrt{1-p^2}$ within the intersection $k$-plane $E_k(p)$, where $p\in[-1,1]$ denotes the distance from $E_k(p)\in{\cal E}_k$ to the origin. A direct calculation shows now
\begin{eqnarray}
\nonumber A(B^d,j,k,m) &=& \int_{-1}^1V_j^m(B^d\cap E_k(p))dp = \int_{-1}^1V_j^m(B_{\sqrt{1-p^2}}^k)dp\\
\nonumber &=& V_j^m(B^k)\int_{-1}^1(1-p^2)^{km\over 2}dp = \left[{k\choose j}{\kappa_k\over\kappa_{k-j}}\right]^m\sqrt{\pi}{\Gamma\left(1+{km\over 2}\right)\over\Gamma\left({3\over 2}+{km\over 2}\right)}.
\end{eqnarray}
To the last formula we can apply Legendre's duplication formula for the Gamma function, which directly proves our claim.\hfill $\Box$\\ \\ The proof of Theorem \ref{thm1} rests upon the following results related to Wiener-It\^o chaos expansions.
\begin{proposition}\label{propMOMENTS}
\begin{itemize}
\item [(a)] The functional $V_{j,k}(W)$ has the Wiener-It\^o chaos expansion \begin{equation} V_{j,k}(W)={\Bbb  E}V_{j,k}(W)+I_1(f_{j,k})\label{eqExpansion}\end{equation} with $f_{j,k}(E_k)=V_j(W\cap E_k)$, $E_k\in{\cal E}_k$.
 \item[(b)] All moments of $I_1(f_{j,k})$ exist and satisfy \begin{equation} {\Bbb  E} I_1(f_{j,k})^m=\tau_k\sum_{i=1}^{m-1} {m-1 \choose i} A(W,j,k,i+1){\Bbb  E}I_1(f_{j,k})^{m-1-i}\label{eqRekursionMomente}\end{equation} for $m\geq 2$.
 \item[(c)] The recursion (\ref{eqRekursionMomente}) has the unique solution \begin{equation}{\Bbb  E}I_1(f_{j,k})^m=\sum_{\pi\in\Pi_m}\tau_k^{|\pi|}\prod_{m_l\in\pi}A(W,j,k,m_l),\label{eqMOMENTS}\end{equation}
for $m\geq 1$.
\end{itemize}
\end{proposition}
Part (a) tells us that $V_{j,k}(W)$ has a rather simple Wiener-It\^o chaos expansion because all kernels $T_nS$ vanish for $n\geq 1$. As a consequence, the central moments of $V_{j,k}(W)$ are given by ${\Bbb  E}I_1(f_{j,k})^m$. As already discussed in the introduction, such powers of multiple Wiener-It\^o integrals can be computed by a so-called diagram formula, see \cite{PT, Sur}. But instead of this technically complicated formula we use another approach and show that the moments must satisfy the recursion (\ref{eqRekursionMomente}) and derive a solution for this equation.
\paragraph{\textsc{Proof of Proposition \ref{propMOMENTS}}} 
\paragraph{(a)} By the definition of the difference operator and the additivity of the intrinsic volumes one has for $E_k\in{\cal E}_k$,
\begin{eqnarray}
\nonumber D_{E_k}V_{j,k}(W)(\eta_k)&=& \sum_{E\in\eta_k\cup\left\{E_k\right\}}V_j(W\cap E)-\sum_{E\in\eta}V_j(W\cap E)\\
\nonumber &=& V_j(W\cap E_k).
\end{eqnarray}
Since $D_{E_k}V_{j,k}(W)(\eta_k)$ depends only on $E_k$ and is independent of the Poisson $k$-flat process $\eta_k$, $$D_{E_k^{(1)},\ldots,E_k^{(n)}}^nV_{j,k}(W)(\eta_k)=0$$ for $E_k^{(1)},\ldots,E_k^{(n)}\in{\cal E}_k$ and all $n\geq 2$. Thus $$T_1V_{j,k}(W)(E_k)={\Bbb  E}D_{E_k}V_{j,k}(W)=V_j(W\cap E_k)=:f_{j,k}(E_k)$$ and
$$T_nV_{j,k}(W)(E_k^{(1)},\ldots,E_k^{(n)})=0$$ for $E_k^{(1)},\ldots,E_k^{(n)}\in{\cal E}_k$ and $n\geq 2$. By (\ref{eqCHAOS}), the Wiener-It\^o chaos expansion (\ref{eqExpansion}) follows immediately.
\paragraph{(b)} It is obvious that all moment of $V_{j,k}(W)$ are finite, because $0\leq V_{j,k}(W)\leq c_{W,j} \eta_k(W)$ with a constant $c_{W,j}>0$ depending on the geometry of the window $W$ and on the parameter $j$.
For the proof of the recursion  we start with the claim that
\begin{equation}D_{E_k}(I_1(f_{j,k})^m)=\sum_{i=1}^m{m \choose i}f_{j,k}(E_k)^iI_1(f_{j,k})^{m-i}\label{eqDiffProduct}\end{equation} for $E_k\in{\cal E}_k$ and with $f_{j,k}(E_k)=V_j(W\cap E_k)$.
This is shown by the product formula (\ref{eqProductDifference}) for the difference operator and by induction. For $m=1$, the statement is easy to verify. For this reason and to prepare the induction step, we start with $m=2$ and obtain $$D_{E_k}(I_1(f_{j,k})^2)=2I_1(f_{j,k})D_{E_k}I_1(f_{j,k})+(D_{E_k}I_1(f_{j,k}))^2=2f_{j,k}(E_k)I_1(f_{j,k})+f_{j,k}(E_k)^2$$ and moreover
\begin{eqnarray*}
& & D_{E_k}(I_1(f_{j,k})^{m+1})\\ &=& (D_{E_k}I_1(f_{j,k}))I_1(f_{j,k})^m+\left[I_1(f_{j,k})+D_{E_k}(I_1(f_{j,k}))\right]D_{E_k}(I_1(f_{j,k})^m)\\ &=& f_{j,k}(E_k)I_1(f_{j,k})^m+(I_1(f_{j,k})+f_{j,k}(E_k))\sum_{i=1}^m{m \choose i}f_{j,k}(E_k)^iI_1(f_{j,k})^{m-i}\\ &=& \sum_{i=1}^{m+1}\left[{m \choose i-1}+{m \choose i}\right]f_{j,k}(E_k)^iI_1(f_{j,k})^{m+1-i}\\ &=& \sum_{i=1}^{m+1}{m+1 \choose i}f_{j,k}(E_k)^iI_1(f_{j,k})^{m+1-i}.\end{eqnarray*}
It follows directly from the series expansion (\ref{eqDiffProduct}) that $T_n(I_1(f_{j,k})^{m-1})=0$ for $n>m-1$. Indeed, each application of the difference operator decreases the power of the Wiener-It\^o integral by one. Thus, after an $m$-fold application of $D$, we get $0$ on the right side of the equation. To proceed, observe that by Lemma \ref{lemItoIntegrals} and the Wiener-It\^o chaos expansion of $I_1(f_{j,k})^m$ combined with the expansion (\ref{eqDiffProduct}), we get
\begin{eqnarray*}
& & {\Bbb  E}I_1(f_{j,k})^m \\ &=& {\Bbb  E}[I_1(f_{j,k})^{m-1}I_1(f_{j,k})]={\Bbb  E}\left[\sum_{n=0}^{m-1}\frac{1}{n!}I_n\left({\Bbb  E}[D_{E_k^{(1)},\ldots,E_k^{(n)}}^n(I_1(f_{j,k})^{m-1})]\right)I_1(f_{j,k})\right]\\ &=& {\Bbb  E}[I_1({\Bbb  E}[D_{E_k}(I_1(f_{j,k})^{m-1}]))I_1(f_{j,k})]=\int_{[W]_k}f_{j,k}(E_k){\Bbb  E}[D_{E_k}(I_1(f_{j,k})^{m-1})]\tau_k\Lambda_k(dE_k)\\ &=& \sum_{i=1}^{m-1}{m-1\choose i}\left(\int_{[W]_k}f(E_k)^{i+1}\tau_k\Lambda_k(dE_k)\right){\Bbb E}[I_1(f_{j,k})^{m-i-1}]\\ &=& \tau_k\sum_{i=1}^{m-1}{m-1\choose i}A(W,j,k,i+1){\Bbb E}I_1(f_{j,k})^{m-i-1}.
\end{eqnarray*}
\paragraph{(c)} The identity
$$\tau_k\sum_{i=1}^{m-1} {m-1 \choose i} A(W,j,k,i+1)\sum_{\pi\in\Pi_{m-i-1}}\tau_k^{|\pi|}\prod_{m_l\in\pi}A(W,j,k,m_l)$$ $$=\sum_{\pi\in\Pi_{m}}\tau_k^{|\pi|}\prod_{m_l\in\pi}A(W,j,k,m_l)$$
follows directly by a decomposition of the partitions of $m$ objects in a set which includes the first object and the partitions of the remaining objects. This finally completes the proof of the proposition.\hfill $\Box$
\paragraph{\textsc{Proof of Theorem \ref{thm1}}}
By (\ref{eqMOMENTS}), we have $$\mu_m(V_{j,k}(W))={\Bbb  E}(V_{j,k}(W)-{\Bbb  E}V_{j,k}(W))^m={\Bbb  E}I_1(f_{j,k})^m=\sum_{\pi\in\Pi_m}\tau_k^{|\pi|}\prod_{m_l\in\pi}A(W,j,k,m_l),$$
which proves the first part of (a).\\
The exact expression for the cumulants follows by a straight forward induction from the expression for $\mu_m(S(W))$ and by using the general well known recursion formula \begin{equation}\gamma_m(V_{j,k}(W))=\mu_m(V_{j,k}(W))-\sum_{i=1}^{m-1}{m-1\choose i-1}\gamma_i(V_{j,k}(W))\mu_{m-i}(V_{j,k}(W)).\label{eqRecCum}\end{equation} 
By the relation (\ref{eqRecCum}) between $\mu_m$ and $\gamma_m$ it is obvious that $\gamma_1(V_{j,k}(W))=\mu_1(V_{j,k}(W))=0$ and $\gamma_2(V_{j,k}(W))=\mu_2(V_{j,k}(W))-0=\tau_k A(W,j,k,2).$
Furthermore, for $m\geq 3$ it holds
\begin{eqnarray*}
\gamma_{m+1}(V_{j,k}(W)) &=& \mu_{m+1}(V_{j,k}(W))-\sum_{i=1}^{m}{m\choose i-1}\gamma_i(V_{j,k}(W))\mu_{m+1-i}(V_{j,k}(W))\\
&=& \sum_{\pi\in \Pi_{m+1}}\tau_k^{|\pi|}\prod_{m_l\in\pi}A(W,j,k,m_l)-\sum_{i=2}^{m}{m\choose i-1}\tau_k A(W,j,k,i)\times\ldots\\
& & \ \ \ \ \ \ \ \ \ \ \ \ldots\times\sum_{\pi\in\Pi_{m+1-i}}\tau_k^{|\pi|}\prod_{m_l\in\pi}A(W,j,k,m_l).
\end{eqnarray*}
The second sum includes all partitions of the first sum with two or more sets, such that
$$\gamma_{m+1}(V_{j,k}(W))=\tau_k A(W,j,k,m+1),$$ completing the proof of part (a).\\ The asymptotic expression for the moments can be derived as follows. By Proposition \ref{propHOMOA}, $A(W,j,k,m_l)$ is homogeneous of degree $jm_l+d-k$, thus $\rho$ appears in the expansion for $\mu_m(V_{j,k}(W_\varrho))$ with power $\displaystyle\sum_{l=1}^{|\pi|}[jm_l+d-k]=jm+|\pi|(d-k)$, where $|\pi|$ is the number of elements of $\pi\in\Pi_m$. For even $m$ this is maximal if $|\pi|=m/2$, i.e. if each $m_l$ in the decomposition of $m$ is $2$. For odd $m$, the decomposition of $m$ into a sum $m_1+\ldots+m_{|\pi|}$ must contain a summand $3$, say $m_1=3$, to be maximal. The remaining summands $m_i$, $i\geq 2$, can then chosen to be equal to $2$. The combinatorial factors $(m-1)!!$ and ${m\choose 3}(m-4)!!$, respectively, simply take into account the number of partitions which realize this situation. The asymptotic formula for the cumulants can be seen from $$\gamma_m(V_{j,k}(W_\varrho))=\tau_k A(W_\varrho,j,k,m)=\tau_k\varrho^{jm+d-k}A(W,j,k,m),$$ where we have used Theorem \ref{thm1} (a) together with Proposition \ref{propHOMOA}. The formulas for the normalized intrinsic volumes $V_{j,k}^*(W_\varrho)$ follow upon division by $({\Bbb V}V_{j,k}(W_\varrho))^{m\over 2}$, in which $\varrho$ appears with power $jm+{m(d-k)\over 2}$. \hfill $\Box$ 
\paragraph{\textsc{Proof of Theorem \ref{thm2}}} As already discussed, there are several ways to see the normal convergence of $V_{j,k}^*(W_\varrho)$. The characteristic function of $V_{j,k}^*(W_\rho)$ can easily be seen to converge to $e^{-{t^2/2}}$, the characteristic function of $\cal N$, by using part (c) of Theorem \ref{thm1}. On the other hand, part (c) of Theorem \ref{thm1} also shows that the cumulants of $V_{j,k}^*(W_\varrho)$ tend to those of a standard normal random variable. Another proof can be based on a general estimate for the Wasserstein distance of $V_{j,k}^*(W_\varrho)$ and $\cal N$ from \cite{P}, which is easily applied in our setting since $V_{j,k}^*(W)$ has a simple Wiener-It\^o chaos expansion. The Berry-Esseen-type estimate is a direct consequence from the result in \cite{Lane}, saying that the rate of convergence in Kolmogorov distance is up to a constant given by
\begin{equation}\label{eqn:BoundsCLTs}
\frac{\tau_k}{\left({\Bbb V} V_{i,k}(W_\varrho)\right)^{\frac{3}{2}}}\int_{[W_\varrho]_k}V_i^3(W_\varrho\cap E_k)\Lambda(d E_k)=\frac{\tau_k A(W_\varrho,j,k,3)}{\left(\tau_k A(W_\varrho,j,k,2)\right)^{\frac{3}{2}}},
\end{equation}
which is of order $\tau_k^{-\frac 12}\varrho^{-\frac{d-k}{2}}$ by Proposition \ref{propHOMOA}.\\
 We turn now to the multivariate central limit theorem, part (b) of Theorem \ref{thm2}. By the chaos expansions of $V_{i,k}(W_\varrho)$ and $V_{j,k}(W_\varrho)$ and the covariance formula (\ref{eqn:CovChaos}) with $n=1$, we obtain
$${\Bbb C}ov(V_{i,k}(W_\varrho),V_{j,k}(W_\varrho))=\langle f_{i,k},f_{j,k}\rangle_1=\int_{[W_\rho]_k}V_i(W_\rho\cap E_k)V_j(W_\rho\cap E_k)\tau_k\Lambda_k(dE_k)$$
and normalization and the homogeneity in $\varrho$ leads to the claimed expression for ${\Bbb C}ov(V^*_{i,k}(W_\varrho),V^*_{j,k}(W_\varrho))$.
In Corllaries 3.4 and 4.3 in \cite{PZ}, multivariate central limit theorems for vectors of first-order Wiener-It\^o integrals are presented as a consequence of more general results. The bounds for two different probability metrics, which both imply convergence in distribution, are up to a constant, depending on the covariance matrix of the multivariate Gaussian distribution, given by the sum over integrals of type (\ref{eqn:BoundsCLTs}) for $i=0,\hdots,d$, which also leads to a rate of convergence of order $\tau_k^{-\frac 12}\varrho^{-\frac{d-k}{2}}$. This immediately implies the statement by noting that the covariances ${\Bbb C}ov(V^*_{i,k}(W_\varrho),V^*_{j,k}(W_\varrho))$ do not depend on the scaling parameter $\varrho$. The case $W=B^d$ can handled similar as in the proof of Lemma \ref{lemmaABdn}. \hfill $\Box$ 

\section*{Acknowledgements}

We are obliged to Matthias Reitzner for his helpful comments and remarks. Moreover, we would like to thank G\"unter Last and Mathew Penrose for making available an earlier version of \cite{LP} to our research group.

\end{document}